\documentclass[twocolumn,10pt]{asme2ej}

\usepackage{epsfig} 
\usepackage{amsmath,amssymb, amsxtra, mathrsfs}
\usepackage{xcolor}
\usepackage{multirow}
\usepackage{caption}
\usepackage{subcaption}
\usepackage{url}

\usepackage{pgfplots,tikz}
\pgfplotsset{compat=1.15}
\usepackage{tikz}
\usetikzlibrary{shapes,arrows,spy,positioning,patterns,shapes.geometric}
\usetikzlibrary{external}
\title{Optimization of District Heating Network Parameters in Steady-State Operation}

\author{Sai Krishna K. Hari\thanks{Address all correspondence to this author.}
    \affiliation{
	Staff Scientist, Member of ASME\\
	Applied Mathematics \& Plasma Physics \,\,\,\\
	Los Alamos National Laboratory\\
	Los Alamos, New Mexico 87545\\
    Email: hskkanth@lanl.gov
    }
}

\author{Anatoly Zlotnik
    \affiliation{
	Staff Scientist\\
	\,\,\, Applied Mathematics \& Plasma Physics\\
	Los Alamos National Laboratory\\
	Los Alamos, New Mexico 87545\\
    Email: azlotnik@lanl.gov
    }
}

\author{Shriram Srinivasan
    \affiliation{
	Staff Scientist\\
	Applied Mathematics \& Plasma Physics \,\,\, \\
	Los Alamos National Laboratory\\
	Los Alamos, New Mexico 87545\\
    Email: shrirams@lanl.gov
    }
}

\author{Kaarthik Sundar
     \affiliation{
	Staff Scientist\\
	  Information Systems \& Modeling \\
	Los Alamos National Laboratory\\
	Los Alamos, New Mexico 87545\\
    Email: kaarthik@lanl.gov
    }
}

\author{Mary Ewers
       \affiliation{
	Staff Scientist\\
	Information Systems \& Modeling\\
	Los Alamos National Laboratory\\
	Los Alamos, New Mexico 87545\\
    Email: mewers@lanl.gov
    }
}

\begin{document}

\maketitle

\begin{abstract}
{\it
We examine the modeling, simulation, and optimization of district heating systems, which are widely used for thermal transport using steam or hot water as a carrier.  We propose a generalizable framework to specify network models and scenario parameters, and develop an optimization method for evaluating system states including pressures, fluid flow rates, and temperatures throughout the network.  The network modeling includes pipes, thermal plants, pumps, and passive or controllable loads as system components. We propose basic models for thermodynamic fluid transport and enforce the balance of physical quantities in steady-state flow over co-located outgoing and return networks.  We formulate an optimization problem with steam and hot water as the outgoing and return carriers, as in legacy 20th century systems. The physical laws and engineering limitations are specified for each component type, and the thermal network flow optimization (TNFO) problem is formulated and solved for a realistic test network under several scenarios.
}
\end{abstract}



\section{Introduction}
\noindent There is significant renewed interest in the use of district heating systems for thermal management worldwide.  District heating can enable the utilization of waste heat generated by thermal power plants, leverage economies of scale, and effectively reduce CO$_2$ emissions \cite{hinkelman2022fast}.  Many university and corporate campuses, government facilities, and urban centers in the United States were developed with district heating in the second half of the 20th century.  Most of these legacy systems (almost 97\%) utilize steam as the carrier \cite{llc2018us}.  Newer systems constructed in the past several decades utilize hot water, which has significant efficiency advantages \cite{novitsky2020smarter}.  In addition to heating buildings, thermal systems can be used for industrial purposes such as wastewater treatment, and there is also increasing consideration of their role in integrated energy network management \cite{schwele2020coordination}.

Modeling of thermal systems has been studied historically in order to aid design and management \cite{marconcini1978numerical}, and physical modeling is well-developed and verified using experiments and observations \cite{sartor2017experimental,gabrielaitiene2006dynamic}.  Simulations and optimization models have been developed for various specific networks in steady-state flow \cite{luo2012steam,wang2016method,stevanovic2007efficient} as well as the in the transient regime \cite{lesko2017modeling,jie2015operation,dahash2019comparative,bohm2000transient,chertkov2019thermal}.  The focus of these studies is usually on specific types of systems and/or detailed physical modeling and model validation.  

In this study, we aim to develop a generalizable framework in which network models and scenario parameters can be specified, and to which computational methods can be applied to readily simulate and optimize system variables including pressures, fluid flow rates, and temperatures throughout the network.  We develop definitions and balance laws for junctions, pipes carrying hot water or steam, thermal plants, pumps, and passive or controllable loads as major system components. These are used to model thermal and fluid transport and conservation in steady-state for outgoing and return flows on co-located networks.  The modeling approach developed here makes use of aspects of our previous work on natural gas transmission networks \cite{hari2021operation,srinivasan2022numerical}.

The aim of this study is to develop generalizable operating modes, physical laws, and engineering limitations for the key types of component present in legacy district heating systems.  These are synthesized to formulate the 
thermal network flow optimization (TNFO) problem where plant settings and controllable loads are optimized according to an economic objective that prioritizes the thermal loads.  The TNFO problem is solved subject to constraints that constitute a well-posed steady-state boundary value problem that has a consistent solution for pressures, flows, and temperatures throughout a thermal network where heat loads and plant output are specified.

The rest of the paper is as follows.  In Section \ref{sec:modeling}, we describe physical modeling of thermal transport systems. In Section \ref{sec:formulation}, we present a network modeling approach for the description of such systems.
In Section \ref{sec:formulation}, we formulate the TNFO steady-state thermal network flow optimization problem, and in Section \ref{sec:computational} we present a computational study in which we solve the TNFO on a test network based on a functioning district heating system under several contingency scenarios. Finally, we summarize the discussion in Section \ref{sec:summary}.


\section{Modeling} \label{sec:modeling}
\noindent 
We model the district heating network as a connected and directed graph, $G(V,E)$. where the edges represent network components such as pipes, plants, and loads, and the vertices ($V$) represent the interconnections between the components. These interconnections are also referred to as junctions.

\subsection{Plant}
\noindent 
The plant receives water mainly from the return system at physically realistic values of temperature and pressure with the primary goal of generating steam at the desired flow rate, temperature, and pressure to meet the district's heating requirements. In modeling the plant, we make three key assumptions. First, we assume that the fluid flow rate in the network can be adjusted \textit{solely} at the plant. Second, the outlet fluid pressure of the plant can be regulated to the desired value using compressors or pumps. Third, we assume that the plant, constrained by an energy production rating, provides the necessary thermal energy to heat the inlet water to convert it into steam at the desired temperature. In this study, we do not model the energy required to adjust fluid pressure; instead, our focus is on analyzing the energy needed to heat water to the desired temperature for steam generation.

We model plants as edges in the graph and denote the set of all plants in the network as $\mathcal{E}_{plant}$; typically, this set is a singleton. Let us consider a plant $e \in \mathcal{E}_{plant}$ that receives water at an inlet temperature $T^{in}_e$ K, pressure $p_{in(e)}$ Pa, and a flow rate $f_e$ Kg/s. Suppose that it is required to supply steam at a setting specified by values for outlet temperature $T^{out}_e$ K, flow rate $f_e$ Kg/s, and pressure $p_{out(e)}$. Then, the process of heating water can be divided into three stages. First, the plant raises the water temperature to $373.15$ K (100 \textdegree C). This stage requires $c^{in}_e f_e (373.15 - T^{in}_e )$ MW of power, where $c^{in}_e$ represents the specific heat capacity of water at constant volume. Second, water is transformed into steam at $373.15$ K, which requires $c^L f_e$ MW of power, where $c^L$ is the latent heat of vaporization. Finally, the steam's temperature is increased to $T^{out}_e$ K at the outlet, which requires $c^{out}_e f_e (T^{out}_e - 373.15)$ MW of power, where $C^{out}_e$ is the specific heat capacity of steam at constant pressure.  Now, suppose that the maximum power rating of the plant is $\mathbf{Pow^{max}_e}$ MW.  Then the total power supplied by the plant is constrained by the inequality
\begin{multline}
c^{in}_e f_e (373.15 - T^{in}_e ) + c^L  f_e + c^{out}_e  f_e  (T^{out}_e - 373.15) \\
\leq \mathbf{Pow^{max}_e}.
\label{eq:max-plant-supply}
\end{multline}

\subsection{Loads}
\noindent 
We refer to facilities with radiators and condensers that extract thermal energy from the carrier and reduce its temperature as \emph{loads}. 
When the outgoing carrier is steam, we suppose that if the temperature reduction does not condense the steam, then the condenser converts the steam to hot water. That is, while fluid enters as steam at the load inlet, it exits as water at the load outlet. Here, we model loads as edges that connect at junctions to pipes carrying steam and water at their inlet and outlet respectively. We denote the set of all loads by $\mathcal{E}_l$, where $\mathcal{E}_l \subset \mathcal{E}$. 

Consider a load $e$, where $e \in \mathcal{E}_l$, and suppose that steam is supplied to the load at a rate of $f_e$ kg/s and a temperature of $T^{in}_e$ K at its inlet. Then, as steam passes through the load, its temperature first reduces to 373.15 K ($100$ \textdegree C). This temperature change releases a power of $c^{in}_e f_e (T^{in}_e - 373.15)$, where $c^{in}_e$ is the specific heat capacity of steam at constant pressure (at the load inlet). Then, at 373.15 K, steam condenses to water, releasing $c^L  f_e$ units of power. Finally, water cools down from 373.15 K to $T^{out}_e$ at the load outlet. This final temperature change in the load releases $c^{out}_e  f_e  (373.15 - T^{out}_e)$ units of power, where $c^{out}_e$ is the specific heat capacity of water at constant volume. Therefore, the total power released at or supplied to the load is  $c^{in}_e f_e (T^{in}_e - 373.15) + c^L  f_e + c^{out}_e  f_e  (373.15 - T^{out}_e)$.

Now, suppose that the power required by the buildings at the load is $Q_e$ MW. The actual power supplied to the load may be either less than or greater than $Q_e$. The former occurs in situations where the total power required by all the buildings exceeds the plant capacity. Suppose that the amount of power that the radiator must supply to the building to maintain its temperature cannot be delivered.  We model the shortfall as undelivered power, denoted by $QS_e$. Conversely, situations where the power required by the buildings at a load is lower than that energy released by condensation may arise to maintain feasibility of thermal and carrier flow balance throughout the network.  In this case we suppose that excess energy, denoted by $QE_e$, is released at the load at the condenser. Then, the power supplied to a load can be related to the power required by the buildings at the load as
\begin{align}
    c^{in}_e f_e (T^{in}_e - 373.15) + c^L  f_e + c^{out}_e  f_e  (373.15 - T^{out}_e) \quad \nonumber \\
     = Q_e + QE_e - QS_e, \label{eq:power_supplied} \\
     QE_e \geq 0, \quad QS_e \geq 0. \qquad \qquad \qquad \label{eq:power_supplied_discrepancies}
\end{align}
The undelivered (shortfall) power $QS_e$ and excess power $QE_e$ are effectively slack variables in the TNFO problem, and are constrained to be positive as reflected by the inequalities \eqref{eq:power_supplied_discrepancies}.  To ensure that the carrier enters the load as steam and exits the load as water, we constrain its temperature at the inlet and outlet of the load by
\begin{align}
T^{in}_e \geq 373.15, \quad   T^{out}_e \leq 373.15. 
    \label{ineq:load-temperatures}
\end{align}
While we do not explicitly quantify the pressure lost in loads, we assume that the pressure at the outlet of a load, $p_{out(e)}$ cannot be greater than that at its inlet, $p_{in(e)}$. For load $e$, this is expressed as
\begin{align}
    p_{in(e)} \geq p_{out(e)}.
    \label{ineq:load-pressures}
\end{align}

\subsection{Pipes}
\noindent 
We consider two types of pipes in the network that comprise outgoing and return systems through which the carrier flows to loads and then returns to the plant. The outgoing system pipes carry steam from the plant to the loads, whereas the return system pipes carry the condensed water back from the loads to the plant. 
We model the pipes as edges and denote the sets of outgoing system and return system pipes by $\mathcal{E}_{op}$ and $\mathcal{E}_{rp}$, respectively. The set $\mathcal{E}_{p}=  \mathcal{E}_{op} \cup \mathcal{E}_{rp}$ denotes the set of all pipes. Because of friction and imperfect insulation in the pipes, the pressure and temperature of the carrier decrease along the direction of the flow.

Consider a carrier flowing along pipe $e \in \mathcal{E}_p$, of length $L_e$, at the rate of $f_e$ kg/s. Suppose that the coefficient of thermal loss of the carrier along the pipe is $\gamma_e$ and the temperature external to the pipe is $T_{ext}$. Then, the decrease in the carrier's temperature along the pipe's length is given by 
\begin{equation}
\label{eq:thermal-loss-ode}
    \frac{dT(x)}{dx} =  - \frac{\gamma_e}{c_e f_e} (T(x) - T_{ext}),
\end{equation}
where $T(x)$ denotes the carrier's temperature at $x$ units from the pipe's inlet along its length, and $c_e$ is the specific heat capacity of the carrier \cite{hinkelman2022fast},  \cite{novitsky2020smarter,chertkov2019thermal}.
Integrating equation \eqref{eq:thermal-loss-ode} leads to a relation between the inlet and outlet temperatures of the pipe, $T^{in}_e$ and $T^{out}_e$ respectively, of form
\begin{align}
\label{eq:temperature-drop}
    T^{out}_e = T_{ext} + (T^{in}_e - T_{ext}) \cdot \exp\left(\frac{-L_e  \gamma_e }{c_e  f_e}\right), \quad \forall e \in \mathcal{E}_p.
\end{align}
Examining the thermal dissipation coefficient $\gamma_e$ in more detail, we suppose that its general form can be determined according to the relation
\begin{equation} \label{eq:heatlosscoeff}
    \frac{1}{\gamma} = \frac{1}{\alpha_a \pi d_a} + \frac{1}{2 \pi K_a} ln \frac{d_b}{d_a} + \frac{1}{2 \pi K_b} ln \frac{d_c}{d_b} + \frac{1}{\pi d_c (\alpha_b + \alpha_c)}.
\end{equation}
Here $\alpha_a$, $\alpha_b$, and $\alpha_c$ are coefficient of forced convection between the steam and the inner wall of the pipe, coefficient of natural convection of the external surface of the steam pipeline, and the radiation coefficient of the external surface of the steam pipeline respectively, and $K_a$, $K_b$ are heat conductivity coefficients for steel and of the insulator respectively \cite{marconcini1978numerical}.  The four terms in equation \eqref{eq:heatlosscoeff} represent (i) forced convection between steam and the internal pipe wall, (ii) heat conduction through the pipe wall, the heat insulator and the aluminium lining, and (iii) natural convection and radiation of the external surface of the pipe towards the atmosphere.

To model steam hydraulics, we assume that the carrier behaves as an ideal gas in the operating range of interest and has a constant compressibility factor. The state equation of steam is then given by
\begin{equation} \label{eq:state-eq}
    p = R T \rho,
\end{equation}
where $p$, $T$, $\rho$, and $R$ denote steam pressure, temperature, density, and specific gas constant, respectively.

We consider steady-state flow where the flow velocity is significantly smaller than the speed of sound.  Here we quantify the decrease in steam pressure along a horizontal pipe $e \in \mathcal{E}_{op}$ of length $L_e$, diameter $d_e$, cross-sectional area $A_e$, and friction factor $\lambda^s_e$.  Using the Darcy-Weisbach model of energy dissipation due to turbulent flow \cite{marconcini1978numerical},  the gradient of pressure in the flow direction is described by
\begin{equation} \label{eq:steam-pressure-drop-ode}
    \frac{dp}{dx} = -\frac{\lambda^s_e}{2 A_e^2 d_e} \frac{f_e |f_e|}{\rho}.
\end{equation}
To solve equation \eqref{eq:steam-pressure-drop-ode} for pressure decrease, we account for the equation of state and thermal effects.
Substituting \eqref{eq:state-eq} in \eqref{eq:steam-pressure-drop-ode} leads to
\begin{equation}
    p \frac{dp}{dx} = -\frac{\lambda^s_e R}{2A_e^2 d_e} f_e|f_e| T.
    \label{eq:interm-deriv-1}
\end{equation}
Integrating equation \eqref{eq:thermal-loss-ode} leads to
\begin{equation}
    T(x) = T_{ext} - \frac{c_ef_e}{\gamma_e} \frac{dT}{dx}.
    \label{eq:rewritten-temperature-drop}
\end{equation}
Then, substituting equation \eqref{eq:rewritten-temperature-drop} into equation \eqref{eq:interm-deriv-1} leads to
\begin{equation}
    p \frac{dp}{dx} = -\frac{\lambda^s_e R}{2A_e^2 d_e} f_e|f_e| \cdot \left(T_{ext} - \frac{c_ef_e}{\gamma_e} \frac{dT}{dx}\right).
    \label{eq:interm-deriv-2}
\end{equation}
Based on the assumption that the steam mass flow rate is constant throughout a pipe, the analytical solution to equation \eqref{eq:interm-deriv-2}, obtained by integrating it from $x = 0$ to $x = L_e$, takes the form 
\begin{align}
\label{eq:steam-pressure-drop}
   \!\!\!\!\!\!\! p_{out(e)}^2 \!-\! p_{in(e)}^2 \!=\! 
    \frac{-\lambda^s_e R}{A_e^2 d_e} f_e |f_e|  (T_{ext} L_e \!+\! \frac{c_e f_e}{\gamma_e} (T^{in}_e \!-\! T^{out}_e)), \nonumber \\  \forall e \in \mathcal{E}_{op}. \!\!\!
\end{align}
The equation \eqref{eq:steam-pressure-drop} above describes the relation between the inlet and outlet pressures, $p_{in(e)}$ and $p_{out(e)}$ respectively, of the outgoing system pipes where steam is the carrier in terms of the flow rate $f_e$ and the inlet and outlet temperatures $T^{in}_e$ and $T^{out}_e$, respectively.

To characterize the pressure changes that occur throughout the return system, we suppose that water is the carrier.  We suppose that water is incompressible, so that its density remains constant through the change between the inlet and outlet pressures  $p_{in(e)}$ and $p_{out(e)}$, respectively, of hot water flowing through a return system pipe $e \in \mathcal{E}_{rp}$.  Supposing a flow rate $f_e$ kg/s through a return pipe of length $L_e$, diameter $d_e$, and with Darcy-Weisbach friction factor $\lambda^w_e$, the flow equation simplifies to
\begin{align}
\label{eq:water-pressure-drop}
    p_{out(e)} - p_{in(e)} = - f_e |f_e| \frac{\lambda^w_e L_e}{2 A_e^2 d_e \rho_e}, \quad \forall e \in \mathcal{E}_{rp}.
\end{align}

\subsection{Pumps}
\noindent 
Water pumps are present at several locations along the return pipeline network to boost the pressure of condensed water to the sufficient levels to ensure that it continues to flow and return to the plant. We suppose here that pumps provide an \textit{additive boost} in pressure that is adjusted in the TNFO problem. In our formulation, we suppose that each pump in the system is located at the inlet of a return pipe, and provides a pressure boost that is represented as a constant additive pressure term, $\alpha_e$, that augments the pipe inlet pressure. As a result, equation \eqref{eq:water-pressure-drop} generalizes to 
\begin{align}
\label{eq:water-pressure-drop-with-pump}
    p_{out(e)} - p_{in(e)} = \alpha_e - f_e |f_e| \frac{\lambda_e L_e}{2 A_e^2 d_e \rho_e}, \quad \forall e \in \mathcal{E}_{rp},
\end{align}
where we set $\alpha_e \equiv 0$ if there is no pump at the inlet of pipe $e \in \mathcal{E}_{rp}$.

\subsection{Junctions}
\noindent 
We refer to the interconnections of network edge components such as plants, pipes, and loads as junctions, and they are modeled as vertices of the graph. At every junction, fluid arrives and leaves through several components at different temperatures and flow rates. However, we assume that the carrier streams arriving at a junction from various edges mix perfectly at the junction so that the  temperature of all flows leaving the junction is equivalent. Let the fluid temperature and pressure at junction $i \in \mathcal{V}$ be denoted by $T_i$ and $p_i$ respectively. Then, because of physical continuity, the temperature at the inlet of all the components carrying fluid away from the junction are constrained to be the same and equal to $T_i$. By denoting the set of all components carrying fluid to and from junction $i$ by $\mathcal{E}_{to}(i)$ and $\mathcal{E}_{fr}(i)$, respectively, the continuity of temperature is constrained as
\begin{align}
\label{eq:mixing}
    T^{in}_e = T_i, \quad e \in \mathcal{E}_{fr}(i).
\end{align}
Moreover, we assume that additional fluid, aside from that flowing through the components, can neither be injected into nor withdrawn from the network at the junctions. TIt follows that mass flow and thermal energy must be conserved at every junction. 
The conservation of mass at junction $i$ can be expressed as 
\begin{align}
\label{eq:mass-conservation}
    \sum_{e \in E_{to}(i)} f_e =  \sum_{\hat{e} \in E_{fr}(i)} f_{\hat{e}}, \quad \forall i\in \mathcal{V},
\end{align}
and the conservation of thermal energy can be stated as
\begin{align}
\sum_{e \in E_{to}(i)} f_e  c^{out}_e  T^{out}_e
= \sum_{\hat{e} \in E_{fr}(i)} f_{\hat{e}}  c^{in}_{\hat{e}}  T^{in}_{\hat{e}}, \quad \forall i\in \mathcal{V}.
\label{eq:thermal-balance}
\end{align}

\section{Thermal Network Flow Optimization}
\label{sec:formulation}
\noindent 
We now formulate a physics-constrained optimization problem to determine controllable parameters in a thermal network that result in a feasible configuration of the physical state across the network. First, we summarize the nomenclature that we use to formulate the problem based on the modeling in the previous section. We then describe the operating constraints that are comprised of inequalities \\

\noindent \textbf{Sets and Functions} \\
$\mathcal{V}$ \quad Set of all nodes/vertices \\
$\mathcal{E}$ \quad Set of all edges \\
$\mathcal{E}_{op}$ \quad Set of all pipes in the outgoing network \\
$\mathcal{E}_{rp}$ \quad Set of all pipes in the returning network \\
$\mathcal{E}_p$ \quad Set of all pipes; $\mathcal{E}_{p} = \mathcal{E}_{op} \cup \mathcal{E}_{rp}$ \\
$\mathcal{E}_l$ \quad Set of all loads \\
$\mathcal{E}_{plant}$ \quad Set of plants in the system (usually a singleton set) \\
$P_{to}(i)$ \quad Set of pipes delivering flow into node $i$, $i \in \mathcal{V}$  \\
$P_{fr}(i)$ \quad Set of pipes taking flow from node $i$, $i \in \mathcal{V}$ \\
$L_{to}(i)$ \quad Set of loads with node $i$, $i \in \mathcal{V}$, as their outlet \\
$PL_{to}(i)$ \quad Set of plants with node $i$, $i \in \mathcal{V}$, as their outlet \\
$L_{fr}(i)$ \quad Set of loads with node $i$, $i \in \mathcal{V}$, as their inlet \\
$PL_{fr}(i)$ \quad Set of plants with node $i$, $i \in \mathcal{V}$, as their inlet \\
$\mathcal{E}_{fr}(i)$ \quad Set of all edges with node $i$, $i \in \mathcal{V}$, as their inlet  \\ 
$\mathcal{E}_{to}(i)$ \quad Set of all edges with node $i$, $i \in \mathcal{V}$, as their outlet \\ 
$in(e)$ \quad Inlet node of edge $e$, $e \in \mathcal{E}$ \\
$out(e)$ \quad Outlet node of pipe $e$, $e \in \mathcal{E}$ \\

\noindent \textbf{Parameters} \\
$A_e$ \quad  Inner cross-sectional area of pipe $e \in \mathcal{E}_p$\\
$c_e$ \quad  Specific heat capacity of carrier in pipe $e \in \mathcal{E}_p$\\
$c^L$ \quad  Specific latent heat of vaporization of water\\
$c^{in}_e$ \quad  Specific heat capacity of carrier at edge $e \in \mathcal{E}$ inlet\\
$c^{out}_e$ \quad  Specific heat capacity of carrier at edge $e \in \mathcal{E}$ outlet\\
$T_{ext}$ \quad  External ambient temperature\\
$\lambda_e$ \quad  Friction factor for fluid flow along pipe $e \in \mathcal{E}_p$\\
$d_e$ \quad  Diameter of pipe $e \in \mathcal{E}_{p}$\\
$L_e$ \quad  Length of pipe $e \in \mathcal{E}_{p}$\\
$\gamma_e$ \quad  Coefficient of thermal heat loss from pipe  $e \in \mathcal{E}_p$\\
$\rho_e$ \quad  Density of carrier in pipe $e \in \mathcal{E}_p$\\
$\lambda_e$ \quad  Darcy-Weisbach resistance for pipe  $e \in \mathcal{E}_{p}$\\
$\mathbf{Pow}^{max}_e$ \quad  Maximum power rating of plant $e \in \mathcal{E}_{plant}$\\
$\bar{\mathbf{T}}$ \quad  Maximum plant outlet temperature\\
$\underline{\mathbf{T}}$ \quad  Lower bound on temperature in the network\\
$\overline{\mathbf{p}}_i$ \quad  Maximum pressure at node $i \in \mathcal{V}$\\
$\underline{\mathbf{p}}_i$ \quad  Minimum pressure at node $i \in \mathcal{V}$ \\
$\overline{\mathbf{\alpha}}_e$ \quad  Maximum pressure boost at inlet of return pipe $e \in \mathcal{E}_{rp}$\\
$R^s$ \quad  Specific gas constant of steam\\
$c^s$ \quad  Specific heat capacity of steam at constant pressure\\
$c^w$ \quad  Specific heat capacity of water at constant volume\\
$\rho^s$ \quad  Density of steam\\
$\rho^w$ \quad  Density of water\\
$Q_e$ \quad  Thermal energy flow demand at load $e \in \mathcal{E}_l$ \\

\noindent \textbf{Variables} \\
$p_i$ \quad  Pressure at node $i \in \mathcal{V}$\\
$T_i$ \quad  Temperature at node $i \in \mathcal{V}$\\
$T^{in}_e$ \quad  Inlet temperature of edge $e \in \mathcal{E}$\\
$T^{out}_e$ \quad  Outlet temperature of edge $e \in \mathcal{E}$\\
$\phi_e$ \quad  Mass flux in pipe $e \in \mathcal{E}_p$\\
$f_e$ \quad  Mass flow rate of carrier in edge $e \in \mathcal{E}$\\
$\alpha_e$ \quad  Pressure boost of pump at inlet of pipe $e \in \mathcal{E}_{rp}$ \\
$QS_e$ \quad Unmet power need at load $e \in \mathcal{E}_l$\\
$QE_e$ \quad  Excess power supplied at load $e \in \mathcal{E}_l$\\

Our aim is to determine plant operating parameters that minimize the deviation between the power required at the loads and the power supplied to them. In order to prevent degeneracy in the optimization solution, i.e., where multiple sets of variable values result in the same objective function value, we augment the objective function with plant outlet temperature, pressure, and flow-rates. Therefore, the objective function of interest is stated as
\begin{equation}
\sum_{e \in \mathcal{E}_{l}} (QE_e + QS_e) +  \sum_{e \in \mathcal{E}_{plant}} (p_{out(e)} + T^{out}_e + f_e + p_{in(e)}  ). 
    \label{eq:objective-function}
\end{equation}
In addition to the constraints that characterize the network flow physics as described in Section \ref{sec:modeling}, we include a collection of engineering constraints that delimit feasible operating regimes for the network.  These are of the form
\begin{align}
    \mathbf{\underline{T}} \leq T_i \leq \bar{\mathbf{T}}, \quad \forall i \in \mathcal{V}, \label{bounds:junction-temperature}\\
    \underline{\mathbf{p}} \leq p_i \leq \bar{\mathbf{p}}, \quad \forall i \in \mathcal{V}, \label{bounds:junction-pressure}\\
    f_e \geq 0, \quad \forall e \in \mathcal{E}, \label{bounds:flow}\\
    T^{out}_e \leq \bar{\mathbf{T}}, \quad \forall e \in \mathcal{E}, \label{bounds:max-temp}\\
    0 \leq \alpha_e \leq \bar{\mathbf{\alpha_e}}, \quad \forall e \in \mathcal{E}_{rp}, \!\!\!\!\!\! \label{bounds:alpha}\\
    QE_e, QS_e \geq 0, \quad \forall e \in \mathcal{E}_l, \!\! \label{bounds:load-power-deviation}
\end{align}
where the inequalities \eqref{bounds:flow} and \eqref{bounds:junction-pressure} represent the operational bounds on the junction temperatures and pressures, respectively, inequality \eqref{bounds:alpha} enforces the non-reversal of flow along edges, and equation \eqref{bounds:max-temp} bounds the maximum temperature in the network from above by enforcing a bound on the plant outlet temperature.  Equation \eqref{bounds:alpha} constrains the pressure boost to be positive in the flow direction and restricts the maximum boost provided by pumps, and \eqref{bounds:load-power-deviation} ensures the non-negativity of excess power dissipated through condensers and undelivered power, respectively.

Collecting the above constraints, the thermal network flow optimization (TNFO) problem is formulated as
\begin{equation} \label{eq:tnfo}
    \begin{aligned}
        & \text{Min} & & \begin{tabular}{@{}c@{}}Sum of plant outlet properties \\ and load power discrepancies: \end{tabular} \textnormal{ \,\, Eqn.\ } \eqref{eq:objective-function},\\
        & \text{s.t.} & & \textnormal{Plant constraint: } \eqref{eq:max-plant-supply}, \\
        & & & \textnormal{Load constraints: } \eqref{eq:power_supplied}, \eqref{ineq:load-temperatures}, \eqref{ineq:load-pressures},\\
        & & & \textnormal{Pipe \& Pump constraints: } \eqref{eq:temperature-drop}, \eqref{eq:steam-pressure-drop}, \eqref{eq:water-pressure-drop-with-pump},\\
        & & & \textnormal{Junction constraints: } \eqref{eq:mixing}, \eqref{eq:mass-conservation}, \eqref{eq:thermal-balance},\\
        & & & \textnormal{Operational bounds: } \eqref{bounds:junction-temperature} - \eqref{bounds:load-power-deviation}.\\
    \end{aligned}
\end{equation}

\noindent In the remainder of the paper, we describe a case study in which we solve the TNFO problem \eqref{eq:tnfo} for a network developed based on a functioning district heating system and examine its condition in a baseline setting as well as for several scenarios.

\section{Computational Study}
\label{sec:computational}
\noindent 
Here we demonstrate the use of the TNFO problem formulation \eqref{eq:tnfo} for approximating optimal plant setpoints for operating a district heating network, and for analyzing the corresponding locational kinematic and thermodynamic properties of the carrier. The formulation \eqref{eq:tnfo} is implemented using JuMP v1.0.1 \cite{Lubin2023}, a package for mathematical optimization embedded in Julia, and solved for a test network using IPOPT \cite{wachter2006implementation}, a primal-dual interior point solver for nonlinear optimization. We developed a network-model for this case study based on the network topology and load profiles for a functioning district heating system. 

\subsection{Network Model}
\begin{figure}
    \centering
    \includegraphics[width = \linewidth]{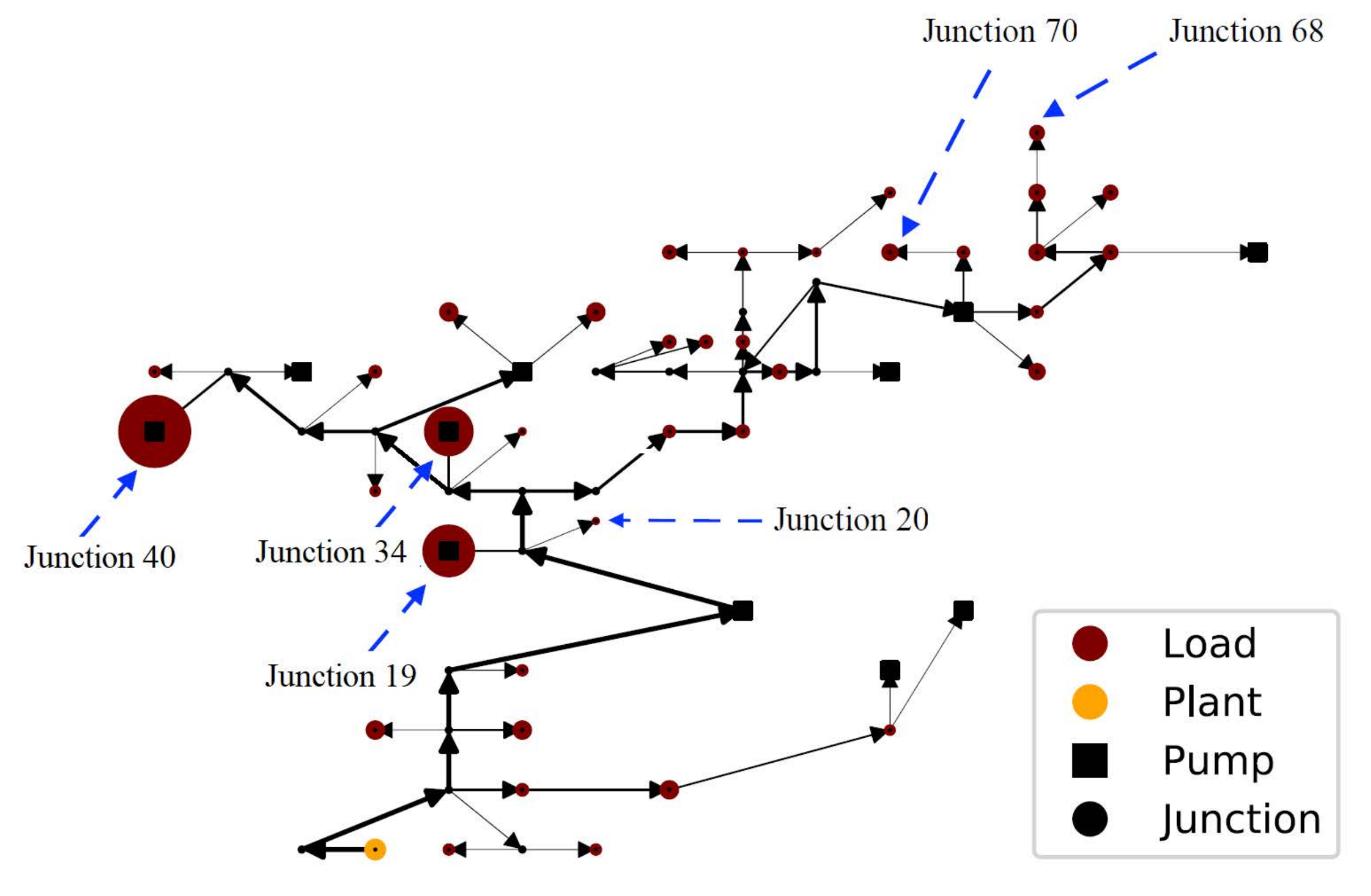}
    \caption{An illustration of the test network used for the computational study. Here, the yellow circle represents the steam plant, the maroon circles represent loads, the black squares represent pumps, and the black edges represent co-located pipes that carry outgoing steam and return water between the plant and the loads. The thickness of each edge is proportional to the diameter of the outgoing system pipe on that edge, and the diameter of each maroon circle is proportional to the amount of thermal energy flow required by the load it represents in the baseline scenario.}
    \label{fig:model}
\end{figure}

\noindent 
The network used for our case study is shown in Figure \ref{fig:model} and consists of the following components:  a steam plant;  $45$ loads;  an outgoing pipeline sub-network comprised of $68$ pipes that carry steam from the plant to various loads;  a returning pipeline sub-network comprised of $68$ pipes that carry condensed water from the loads back to the plant; $11$ pumps located at the inlet of different return pipes; and $134$ junctions that represent the intersections of these components. 
The lengths of the edges in Figure \ref{fig:model} are not drawn to scale.  The thicknesses of the edges illustrated in Figure \ref{fig:model} denote the relative diameters of the pipes in the outgoing network.  The junctions that we refer to in the subsequent discussion are labeled in the figure.

\subsection{Baseline Scenario}
\noindent 
We define a baseline scenario of parameters, thermal loads, and operating setpoints.  We assume that, under normal operational conditions, the plant can provide a maximum power of $30$ MW. In order to determine suitable parameters for the network, we consider a representative day during which the buildings at the loads  require a total of $15.14$ MW of power. Combining this estimate with the relative sizes of the buildings, we construct a \textit{baseline scenario}, where the relative power requirements at the loads are represented by the size of the nodes in Figure \ref{fig:model}.  The three largest circles correspond to three large buildings served by the district heating system.

Next, we assume carrier properties as follows.  Suppose that $R^s = 461.5$ J/kg$\cdot$K, $c^s =1996$ J/kg$\cdot$K, $c^w = 4186$ J/kg$\cdot$K, $c^L = 2230$ kJ/kg, $\rho^s = 0.5$ kg/m$^3$, $\rho^w = 1000$ kg/m$^3$. In order to obtain realistic temperature and pressure profiles for the baseline scenario, we set $\lambda^s_e = 0.01$ and $\gamma_e = 0.1$ for the outgoing system pipes and  $\lambda^w_e = 0.002$ and $\gamma_e = 0.05$ for the return system pipes. For all the scenarios considered in our study, we assume the network operating bound limits as follows.  We suppose that $\bar{\mathbf{T}} = 150 $ \textdegree C, $\mathbf{\underline{T}} = 80$ \textdegree C, $T_{ext} = 25 $ \textdegree C, $\bar{\mathbf{p}} = 80$ psi, and $\mathbf{\underline{p}} = 5$ psi, and we require a minimum plant inlet pressure of $40$ psi. We allow a maximum pressure boost of $5$ psi at the pumps.
With these parameter values, we solve the TNFO problem \eqref{eq:tnfo} to determine the optimal operating parameters and the associated pressure, temperature, and flow solution throughout the network for the baseline scenario. 

\begin{figure*}[h!]
\centering
\begin{subfigure}{0.48\textwidth}
\centering
    \includegraphics[clip, trim=.9cm .9cm .9cm .9cm, width = \linewidth]{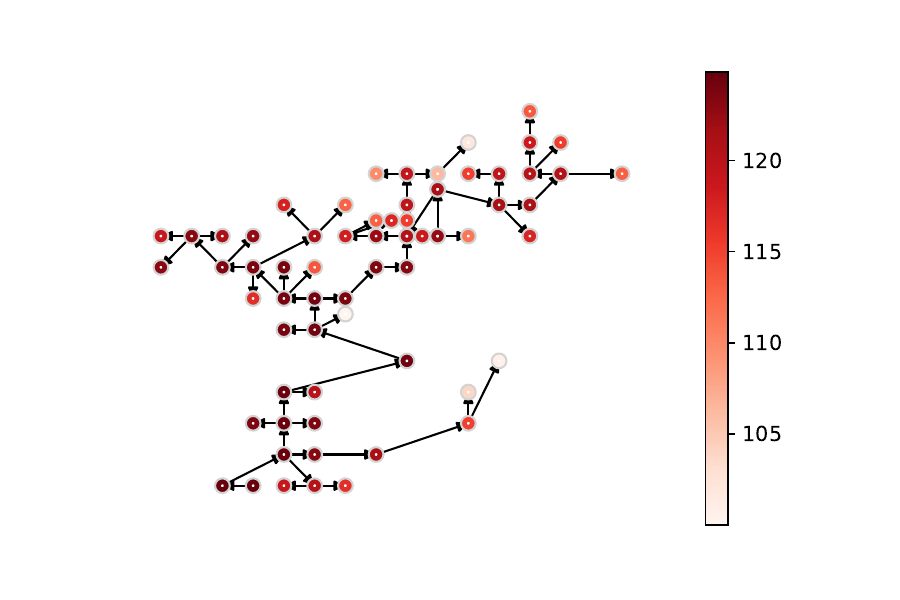}
    \caption{Outgoing temperatures}
    \label{fig:baseline-ot}
\end{subfigure}
\hfill
\begin{subfigure}{0.48\textwidth}
\centering
    \includegraphics[clip, trim=.9cm .9cm .9cm .9cm, width = \linewidth]{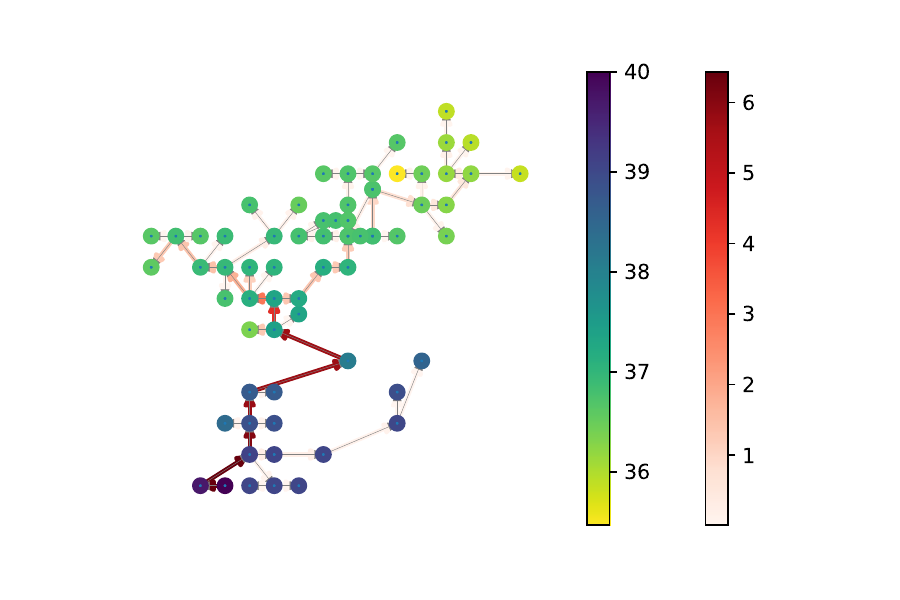}
    \caption{Outgoing pressures and flows}
    \label{fig:baseline-outflow}
\end{subfigure}
\hfill
\begin{subfigure}{0.48\textwidth}
\centering
    \includegraphics[clip, trim=.9cm .9cm .9cm .9cm, width = \linewidth]{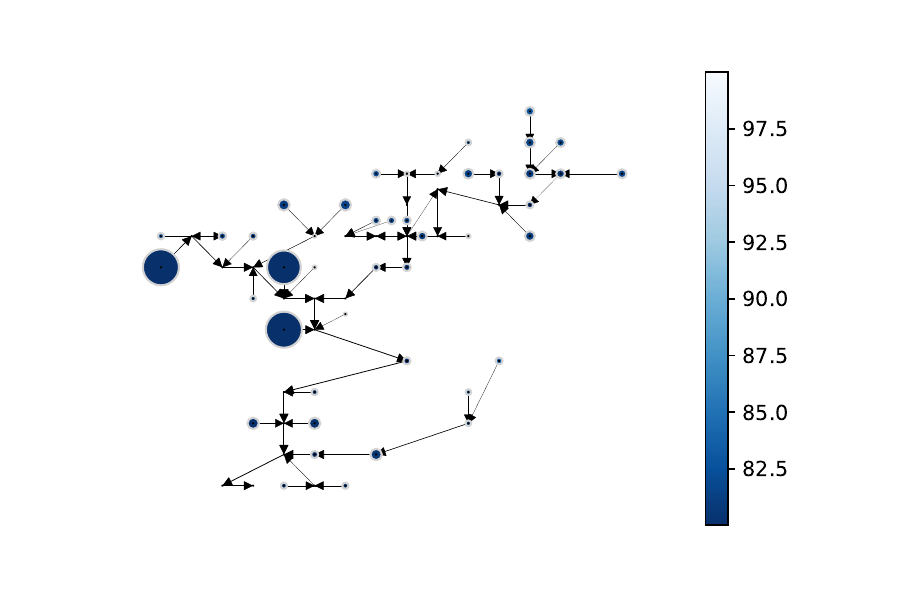}
    \caption{Returning temperatures}
    \label{fig:baseline-rt}
\end{subfigure}
\hfill
\begin{subfigure}{0.48\textwidth}
\centering
    \includegraphics[clip, trim=.9cm .9cm .9cm .9cm, width = \linewidth]{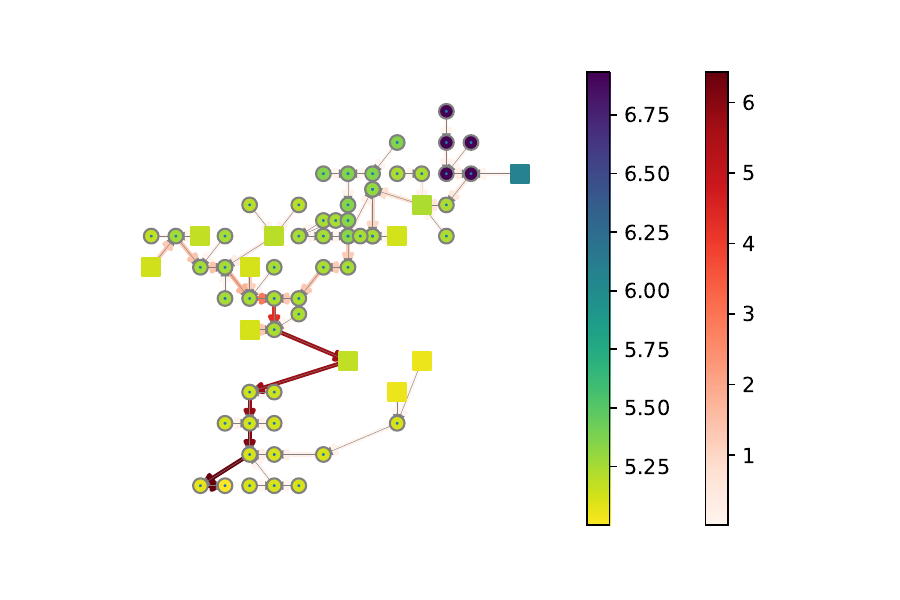}
    \caption{Returning pressures}
    \label{fig:baseline-rp}
\end{subfigure}  
\caption{Baseline Scenario: Optimal physical temperature, pressure, and flow solution.}
\label{fig:baseline}
\end{figure*}

The optimal solution for plant operating parameters given the baseline loads  are output temperature $T_{plant}^{out}= 124.86$ \textdegree C, output pressure $p_{out(plant)}=40$ psi, and mass flow rate of $f_{plant}=6.43$ kg/s. These operating parameters require 15.2 MW of power that is applied to convert condensed water entering the plant inlet to steam of the desired properties leaving the plant outlet, and a total of 0.06 MW of this power is dissipated in pipes. We present the properties of steam and water across the network corresponding to this solution in Figure \ref{fig:baseline}.
The results in Figure \ref{fig:baseline} illustrate the physical temperature, pressure, and flow configuration obtained in the optimal solution to problem \eqref{eq:tnfo} for the baseline scenario.  

Figure \ref{fig:baseline-ot} shows the temperature profile of steam throughout the outgoing network. A red color gradient is used such that the dilution of color represents a decrease in temperature. In the baseline scenario, the temperatures range from 124.86 \textdegree C at the plant outlet, shown in dark red, to 100 \textdegree C at junction 20 (refer to Fig. \ref{fig:model} for junction labels).
In Figure \ref{fig:baseline-outflow}, we use a similar color scheme for edges to illustrate the steam flow rate along the outgoing pipes. Here, darker red corresponds to the flow rate of 6.43 kg/s at the plant outlet and lighter shades represent lower flow rates through branches that distribute the carrier to the loads. In this figure, we also depict the steam pressure at various junctions of the outgoing network. The pressure here ranges from 40 psi at the plant to 35.47 psi at junction 70, and we represent this change using the viridis color gradient, with the colors transitioning from dark purple to green to yellow with decreasing pressure.
The decrease in steam temperature and pressure along the direction of flow is apparent in Figures \ref{fig:baseline-ot} and \ref{fig:baseline-outflow}, and this is consistent with the dissipative effects enforced by constraints \eqref{eq:temperature-drop}, \eqref{eq:steam-pressure-drop}, and \eqref{eq:water-pressure-drop}.

In Figure \ref{fig:baseline-rt}, we illustrate the temperature of hot water in the return system resulting from condensation of steam at the load outlets. We use a blue color gradient to depict the temperatures, with darker blue corresponding to colder temperatures. While the temperatures here range from $100$ \textdegree C to 80.04 \textdegree C, temperature does not monotonically decrease along paths through the return network, as it does along paths through the outgoing network. This can be attributed to the variation in power consumption at different loads and the resulting differences in the load outlet carrier temperatures. The lowest temperature of 80.04 \textdegree C is observed at the outlet of the load at junction 40 (refer to Fig. \ref{fig:model} for junction labels), which is the largest load in the baseline scenario. The temperatures at the other two major loads at junctions 24 and 19 are within 0.04 \textdegree C of that at junction 40. Nevertheless, after the release of hot water at different temperatures, the solution shows perfect mixing at the junctions of the return network as enforced by the model constraints, and the temperatures subsequently decrease along the direction of flow.

In Figure \ref{fig:baseline-rp}, we illustrate the flow rate and pressure of hot water throughout the return network using the same color scheme as in Figure \ref{fig:baseline-outflow}. The flow direction of hot water is opposite to that of steam, and as water flowing towards the plant inlet merges through different branches, the plant through-flow rate of 6.43 kg/s is accumulated. This flow direction reversal appears as a reversal of the sign of the pressure gradient along the return network, which is depicted by the lighter color at the junctions from dark purple to green to yellow. A pressure boost is observed at a few locations due to the presence of pumps, which are indicated by squares at the junctions.
The return network pressures are bounded above by the pressures in the outgoing network by the constraint \eqref{ineq:load-pressures} and bounded below at 5 psi by constraint \eqref{bounds:junction-pressure}. These bounds, together with the pressure decrease along pipes, pressure boost at the pumps, and the minimization of plant inlet pressure, cause the pressures to range from 6.93 psi at junction 68 to 5 psi at the plant inlet. Throughout the rest of the paper, we will use the above results as benchmarks.

\subsection{Contingency Analysis Scenarios}
\noindent 
We now examine how the solutions of the TNFO problem \eqref{eq:tnfo} change in scenarios that represent contingencies such as large increases in load or equipment outages at the plant.  This demonstrates the utility and adaptability of the mathematical programming formulation to examine various operating conditions.   The same network model as illustrated in Figure \ref{fig:model} is subjected to four operational scenarios, which are represented using the following changes to parameter values:
\begin{enumerate}
    \item \textbf{Functional contingency}, which encompasses contingencies such as leaving the door of a large building open: $200 \%$ increase in the thermal load at junction 40, which is the largest load in the network. 
    \item \textbf{Extreme load}, which represents an increase in thermal load that may occur during extreme weather: $50 \%$ increase in the power required at all the loads. 
    \item \textbf{Functional contingency during extreme load}, which is a combination of the above two scenarios: $200 \%$ increase in the load at junction 40 and $50 \%$ increase in the other loads.
    \item \textbf{Equipment outage}, such as maintenance of one of the boilers at the plant: Maximum power supply curtailed to $20$ MW.
\end{enumerate}

For each of the above scenarios, the TNFO program \eqref{eq:tnfo} is solved and the optimal plant operating parameters are tabulated in Table \ref{tab:summary_of_scenarios}. In the rest of the section, we analyze these solutions and explain the deviations in various flow properties that arise because of the plant parameter changes.

\begin{table*}[h!]
\caption{Summary of TNFO problem solutions for baseline and contingency scenarios.}
\vspace{-3ex}
\begin{center}
\label{tab:summary_of_scenarios}
\begin{tabular}{|c|c|cccc|ccc|}
\hline
& Scenario & \multicolumn{4}{c|}{Power} & \multicolumn{3}{c|}{Plant} \\ \cline{3-9} & &Required & Supplied  &  Pipe Losses & Unmet Load  & Temperature  & Pressure  & Flow \\
 & & (MW) & (MW) &  (MW) & \% & (\textdegree C) &  (psi) & (kg/s) \\ \hline
 1 & Baseline & 15.14 & 15.2  & 0.06 & 0        & 124.86 & 40.0 & 6.43  \\ \hline
 2 & junction 40 $\nearrow$200\%     & 21.43 & 21.49  & 0.06 & 0        & 124.66 & 40.0 & 9.10  \\ \hline
3 & \qquad All loads $\nearrow$50\%  & 22.70 & 22.76 & 0.06 & 0        & 115.77 & 40.0 & 9.71  \\ \hline 4 &
\begin{tabular}[c]{@{}c@{}}   junction 40 $\nearrow$200\%,\\  other loads $\nearrow$50\% \end{tabular} & 27.43 & 27.48 & 0.06 & 0        & 115.69 & 40.0 & 11.72 \\ \hline
5 &
\begin{tabular}[c]{@{}c@{}}  Supply curtailed, \\  junction 40 $\nearrow$200\%,\\  other loads $\nearrow$50\% \end{tabular} & 27.43 & 20    & 0.06 &  27.28\% & 115.81 & 40.0 & 8.53  \\ \hline
\end{tabular}
\end{center}
\end{table*}

\subsubsection{Scenario 1: Functional Contingency}
\noindent 
The functional contingency scenario represents situations such as leaving the door of a large structure open.  For the present case study, we model a functional contingency as a $200 \%$ increase in the thermal load at junction 40, which is the largest load in the network.  The optimal solution indicates an increase in plant operating parameters to a flow rate of 9.10 kg/s, at the same temperature and pressure at the plant outlet as the baseline scenario. This increases the plant power output to 21.43 MW from 15.14 MW in the baseline scenario (the additional 6.29 MW is a 40\% increase). The thermal energy dissipated from the outgoing and return system pipes is the same as in the baseline scenario.

\begin{figure}[th!]
\centering
    \includegraphics[clip, trim=1.5cm .9cm .9cm .9cm, width = 0.5\textwidth]{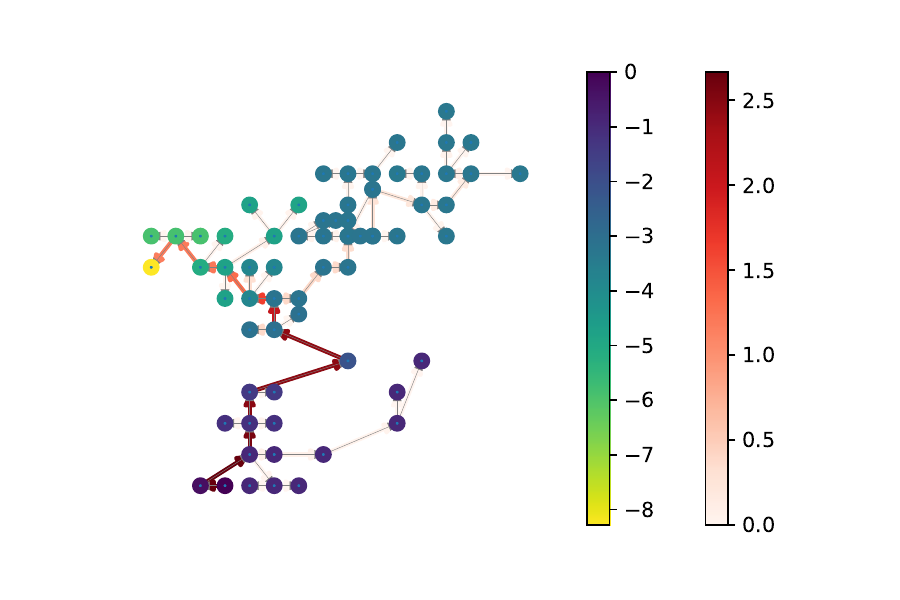}
    \caption{Scenario 1: Change in steam pressures (at junctions) and flow rates (at edges) relative to those in the baseline scenario.}
    \vspace{-2ex}
    \label{fig:case-1-relative-flow}
\end{figure}

The majority of the increased steam flow in this scenario is directed toward junction 40. This can be observed by the red-colored edges in Figure \ref{fig:case-1-relative-flow}, which depicts the change in steam flow rates and pressures in the outgoing system pipes relative to the baseline scenario. Inspection of the figure also suggests that this increased flow rate causes a greater decrease in the network pressures. This is particularly evident for junction 40, at which the steam pressure decreases by 8.29 psi with respect to the baseline scenario (from 36.6 psi to 28.31 psi). Consequently, the junction with the minimum outgoing system pressure also switches from junction 70 (farthest from the plant) in the baseline scenario to junction 40 in this scenario.

The increased steam flow rate corresponds to a greater return water flow rate, which results in a more significant pressure decrease in the return network. However, given that the minimum return system pressure in the baseline scenario attains its lower bound of 5 psi, the pressure decrease in this scenario is accommodated by increasing the maximum return network pressure to 7.16 psi (from 6.93 psi in the baseline scenario). Because the plant output temperature is unchanged in the TNFO solution for this scenario, the average outgoing temperatures remain close to their baseline values. Still, the increased heating need at junction 40 requires the difference between plant inlet and outlet temperatures to increase somewhat compared with the baseline scenario.

\begin{figure}[t!]
    \centering
    \vspace{-2ex}
    \includegraphics[clip, trim=1.5cm .9cm .9cm .9cm, width = 0.5\textwidth]{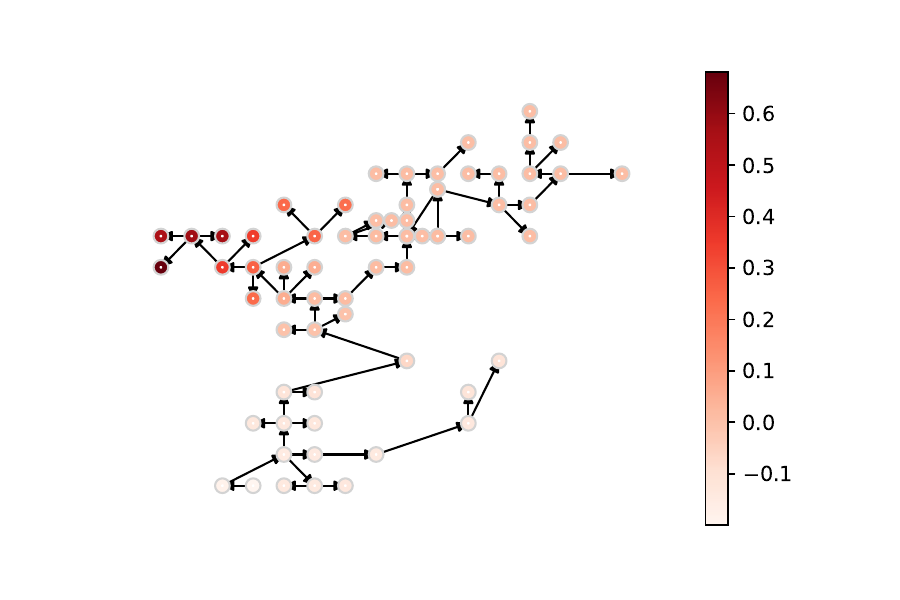}
    \caption{Scenario 1: Increase in the outgoing temperatures at junction 40 with respect to the baseline scenario.}
    \label{fig:case-1-relative-ot}
    \vspace{-2ex}
\end{figure}

\subsubsection{Scenario 2: Extreme Load}
\begin{figure}[h!]
    \centering
    \includegraphics[width = 0.5\textwidth]{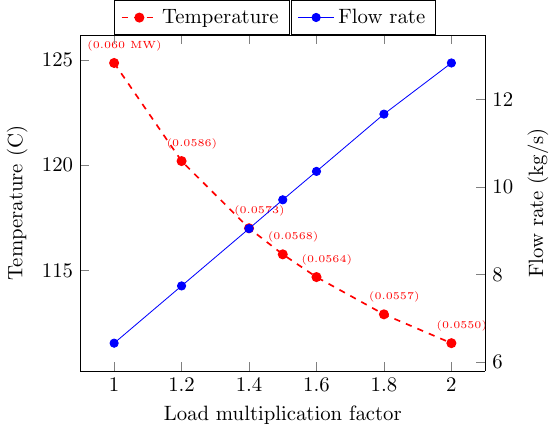}
    \caption{Change in optimal plant operating parameters resulting from a uniform network-wide increase in thermal load. Red and blue lines indicate the plant output temperature and through-flow, respectively, and data point labels in parenthesis indicate the thermal energy in MW dissipated from pipes.}
    \label{fig:optimal-setpoints-all-loads}
\end{figure}

\noindent 
The extreme load scenario represents an increase in thermal load that may occur during extreme weather.  For our case study, we model an extreme load scenario as a $50 \%$ increase in the thermal power required at all the loads in the network.  The total power required at the loads in this scenario is similar to that in Scenario 1. However, unlike the thermal load increase at a single junction, in this scenario all thermal loads are uniformly increased by 50\% from their baseline values. The TNFO simulation solution for this scenario results in a decrease of the plant outlet temperature while increasing the flow rate. The optimal plant operating parameters are a flow rate of 9.71 kg/s, outlet temperature of 115.77 \textdegree C, and outlet pressure of 40 psi.  The latter is the same outlet pressure as in the baseline scenario solution. While the flow rate here does not significantly deviate from that in the solution for Scenario 1, the outlet temperature is 9.09 \textdegree C lower than in the baseline scenario solution. To further examine the implications of the decreased plant temperature, we perform a sensitivity analysis by varying the uniform network-wide load increase percentage from 0\% to 100\% with respect to the baseline scenario load, expressed as a multiplication factor. The resulting optimal plant temperature and flow parameters for this sensitivity study are plotted in Figure \ref{fig:optimal-setpoints-all-loads}, and the thermal losses from pipes are indicated there as well. The results of the sensitivity analysis indicate that lower plant temperatures with increased flow rates are sufficient to service a uniform network-wide thermal load increase. This trend can partially be attributed to the decrease in thermal losses in pipes (displayed inside parenthesis in the figure) with lower supply temperatures, in addition to the nonlinear and non-convex constraints that relate pressures, flow, and thermal dissipation in the pipes.


\begin{figure}[t!]
    \centering
    \includegraphics[clip, trim=1.5cm .9cm .9cm .9cm, width=0.5\textwidth]{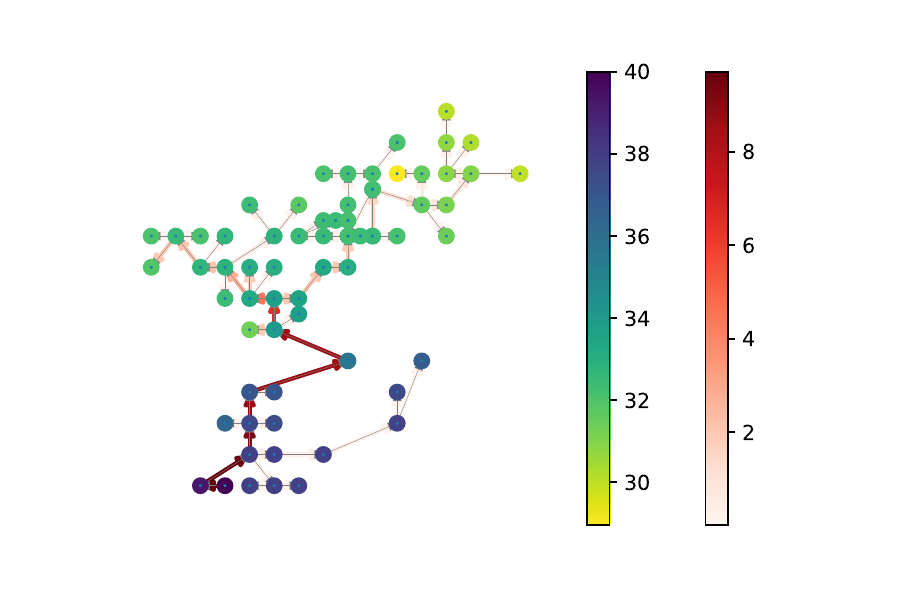}
    \caption{Scenario 2: Locational pressures and flows in the outgoing system.}
    \label{fig:case-2-op}
\end{figure}

Because of the differences in the thermal load distribution between Scenarios 1 and 2, the flow rate patterns are also different. This is illustrated in Figure \ref{fig:case-2-op}, where there is a uniform increase in the flow rate to all the loads, unlike the localized flow increase in Scenario 1. As a result of this uniform flow increase, the pressure decrease pattern in the outgoing and returning networks in this scenario is similar to the baseline scenario. Nevertheless, because of the higher flow rate, the pressure decrease is more significant.
The outgoing pressures range between 40 psi (at the plant outlet) to 28.96 psi (at junction 70), whereas the return pressures range between 9.39 psi (at junction 68) to 5 psi (at the plant inlet).  While the lower plant outlet temperature results in lower outgoing system temperatures, the temperature decrease in the network is similar to that in the baseline scenario.

\subsubsection{Scenario 3: Contingency During Extreme Load}
\noindent This scenario in which a functional contingency occurs during extreme load synthesizes Scenarios 1 and 2, and is modeled for our case study by a $200 \%$ increase in the load at junction 40 and $50 \%$ increase in all other loads.  A combination of the effects observed in the previous two scenarios is reflected in the results.  The total thermal load in this scenario is 27.43 MW, which is the highest among the four scenarios that we consider. To service this load, the optimal solution of the TNFO problem results in plant operating parameters consisting of a flow rate of 11.72 kg/s, outlet temperature of 115.69 \textdegree C, and outlet pressure of 40 psi. Here, the flow rate is the highest and the outlet temperature is the lowest of the solutions obtained in the examined scenarios.

\begin{figure}[t!]
\centering
    \includegraphics[clip, trim=1.5cm .9cm .9cm .9cm, width=0.5\textwidth]{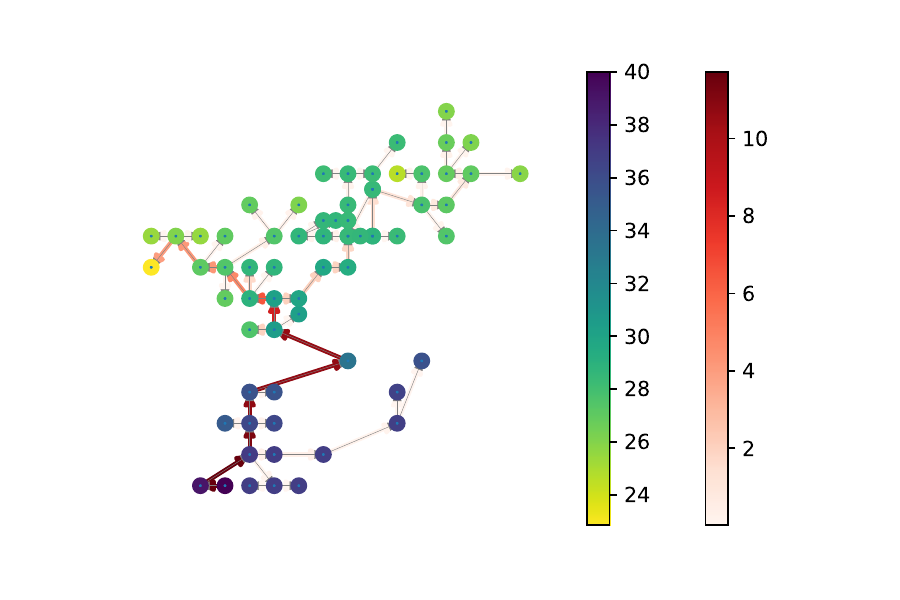}
    \caption{Scenario 3: Locational pressures and flows in the outgoing system.}
    \label{fig:case-3-op}
\end{figure}

The combined effects of the increased loads become apparent, first in the flow rates throughout the network as illustrated in Figure \ref{fig:case-3-op}. Similarly to Scenario 2, the flow rate of the carrier to all loads increases. An even higher increase in the flow is diverted toward junction 40, which has the load with the highest demand, as in Scenario 1. This elevated flow rate is consistent with previous cases and results in greater decreases in pressure throughout the network.

We also observe the effects of increased load in the pressure decreases throught the network. Although junction 40 has the lowest pressure on the outgoing system as in Scenario 1, the next lowest outgoing system pressure occurs at junction 70 where outgoing system pressure is the lowest in Scenario 2. Notably, the optimal plant outlet temperature in this scenario is the lowest among all cases, which in turn leads to lower overall network temperatures with respect to the baseline scenario, similar to Scenario 2. However, the inlet temperature at junction 40 is slightly higher than in Scenario 2, which is caused by the more significant localized increase in load at that junction.  As in previous scenarios, the return network pressures are higher with respect to the baseline scenario solution, with values in the range from 9.63 psi to 5 psi.

\subsubsection{Scenario 4: Equipment Outage}
\begin{figure}[h!]
    \centering
    \includegraphics[clip, trim=1.5cm .9cm .9cm .9cm, width = 0.5\textwidth]{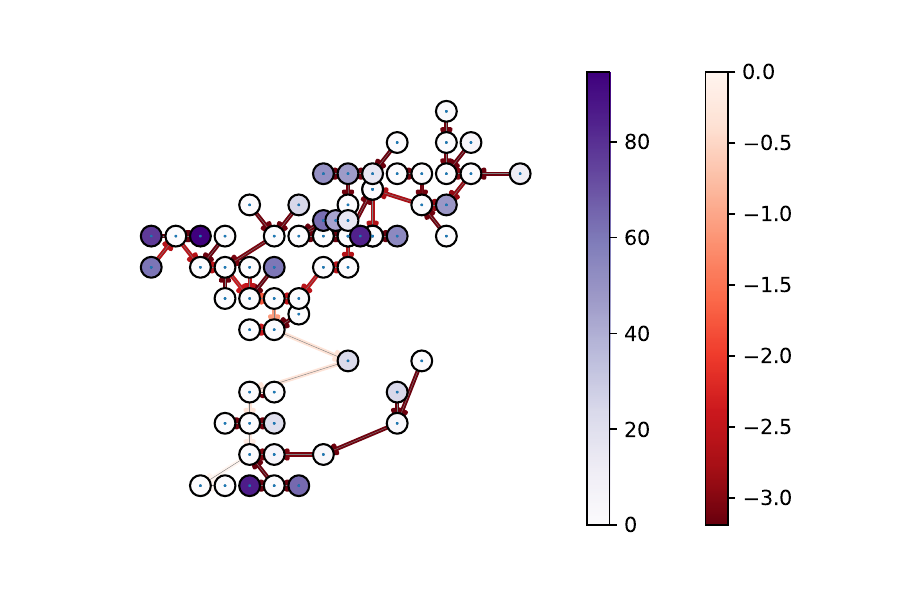}
    \caption{Results of Scenario 4 compared to those of scenario 3. Specifically, the edge color depicts the change in the flow rate and the node color depicts the percent of thermal load demand that is unmet.}
    \label{fig:case-7-ls}
    \vspace{-3ex}
\end{figure}

\noindent 
The final scenario that we consider is an equipment outage, such as in planned or unplanned maintenance of one of the boilers at the plant.  We represent such an event in our case study models as a decrease of the maximum plant power output constraint bound value to $20$ MW.  In this scenario, while the total thermal load in the network is 27.43 MW, the plant supply is limited at 20 MW. As a consequence, the steam flow rate at the plant outlet is reduced to 8.53 kg/s, and thermal energy delivery is insufficient at several loads. The absolute and relative values of unmet thermal energy demand is shown in Figure \ref{fig:case-7-ls}. In this scenario, a total of 27.28\% of thermal demand is unmet, and the distribution of this unmet load throughout the network is determined based on the objective of minimizing this total amount.

\section{Summary}
\label{sec:summary}
\noindent 
We have developed a generalizable framework for optimizing operating parameters for a district heating system that uses steam as the outgoing carrier and hot water as the return carrier.  Given a network model and scenario parameters, our thermal network flow optimization (TNFO) formulation is used to evaluate system states including pressures, fluid flow rates, and temperatures to minimize unmet thermal load throughout the network.  We propose basic models for thermodynamic fluid transport and enforce the balance of physical quantities in steady-state flow over co-located outgoing and return networks.  In addition to first-principles definition of physical laws and specification of engineering limitations in a well-posed optimization formulation, we examine the physical solutions produced by solving the TNFO problem for several scenarios for a realistic network case study, and perform a sensitivity analysis.

Our optimization formulation is developed as an entirely continuous problem, without mixed integer variables for elements such as valves, and with simplified modeling of loads that could be large and complex facilities as thermal sinks and condensers. The solution is driven by the objective function that minimizes unmet demand and excess energy dissipation, and aims to service thermal loads using parameters that minimize plant outputs including outlet pressure and temperature and through flow.  Slack variables for unmet demand and excess energy dissipation ($QS_e$ and $QE_e$, respectively) are included in the formulation to accommodate situations where nominal thermal loads are not met or thermal energy dissipated at a load must be greater than building needs in order to maintain a feasible physical solution.  These slack variables are positive only when required because of network topology limitations.  We observe some counter-intuitive behavior in our sensitivity study where increasing load results in lower plant outlet temperature, but this is explained by increasing carrier flow rates that result in greater thermal energy transfer.  Future studies could examine extensions of our formulation to accommodate integer variables for valves resulting in mixed-integer linear programs, to incorporation of resilience constraints so that the solution delivers sufficient energy to loads given outages or contingencies, to transient flow conditions, and interaction with interconnected infrastructure networks.

\begin{acknowledgment}
This study was supported by the U.S. Department of Energy's Grid Modernization Laboratory Consortium (GMLC) Project on Energy Resilience for Mission Assurance (ERMA). Research conducted at Los Alamos National Laboratory is done under the auspices of the National Nuclear Security Administration of the U.S. Department of Energy under Contract No. 89233218CNA000001. This article has been approved for unlimited release with report no. LA-UR-24-23145.
\end{acknowledgment}

\vspace{2ex}
\bibliographystyle{asmems4}
\bibliography{references}

\end{document}